\documentclass[onefignum,onetabnum]{siamonline250211}

\usepackage{graphicx} 
\usepackage{amsmath}
\newtheorem{conjecture}[theorem]{Conjecture}
\usepackage{amssymb}
\usepackage{url} 
\usepackage{booktabs} 
\usepackage{mathrsfs} 
\usepackage{xcolor} 
\definecolor{qedblue}{HTML}{005188}

\usepackage{enumitem}
\usepackage{algorithm}
\usepackage{algpseudocode}

\usepackage[numbers]{natbib}
\usepackage{geometry}
\usepackage{blindtext}
\usepackage{xr-hyper}

\headers{Theory of Collinearity on Variable Importances}{K. K. Bladen, D. R. Cutler, A. Wisler}

\title{Mathematical Theory of Collinearity Effects on Machine Learning Variable Importance Measures\thanks{Submitted to the editors Sept 30, 2025.}}

\author{Kelvyn K. Bladen\thanks{Department of Mathematics \& Statistics, Utah State University, Logan, UT 
  (\email{kelvyn.bladen@usu.edu}, \email{richard.cutler@usu.edu}, \email{alan.wisler@usu.edu}).} \and D. Richard Cutler\footnotemark[2]
  \and Alan Wisler\footnotemark[2]}

\title{Mathematical Theory of Collinearity Effects on Machine Learning Variable Importance Measures
}

\author{Kelvyn K. Bladen\thanks{Department of Mathematics \& Statistics, Utah State University, Logan, UT 
  (\email{kelvyn.bladen@usu.edu}, \email{richard.cutler@usu.edu}, \email{alan.wisler@usu.edu}).} \and D. Richard Cutler\footnotemark[2]
  \and Alan Wisler\footnotemark[2]}

\begin{document}

\maketitle

\begin{abstract}
In many machine learning problems, understanding variable importance is a central concern. Two common approaches are Permute-and-Predict (PaP), which randomly permutes a feature in a validation set, and Leave-One-Covariate-Out (LOCO), which retrains models after permuting a training feature. Both methods deem a variable important if predictions with the original data substantially outperform those with permutations. In linear regression, empirical studies have linked PaP to regression coefficients and LOCO to $t$-statistics, but a formal theory has been lacking. We derive closed-form expressions for both measures, expressed using square-root transformations. PaP is shown to be proportional to the coefficient and predictor variability: $\text{PaP}_i = \beta_i \sqrt{2\operatorname{Var}(\mathbf{x}^v_i)}$, while LOCO is proportional to the coefficient but dampened by collinearity (captured by $\Delta$): $\text{LOCO}_i = \beta_i (1 -\Delta)\sqrt{1 + c}$. These derivations explain why PaP is largely unaffected by multicollinearity, whereas LOCO is highly sensitive to it. Monte Carlo simulations confirm these findings across varying levels of collinearity. Although derived for linear regression, we also show that these results provide reasonable approximations for models like Random Forests. Overall, this work establishes a theoretical basis for two widely used importance measures, helping analysts understand how they are affected by the true coefficients, dimension, and covariance structure. This work bridges empirical evidence and theory, enhancing the interpretability and application of variable importance measures.
\end{abstract}



\section{Introduction}
\label{s1:intro}

In regression problems, and in supervised learning more broadly, understanding the contribution of different variables to a predictive model is a longstanding challenge. In many cases, assessing the relative importance of variables to a prediction can be as vital to the analysis as the accuracy of the predictions themselves. Consequently, methods for measuring variable importance have been employed across a wide range of disciplines, including social sciences \citep{blalock1961evaluating,darlington1968multiple}, medicine/bioinformatics \citep{Strobl,nicodemus2011stability, goldstein2017opportunities}, ecology \citep{ecology, murray2009methods}, and economics \citep{varian2014big}.

Some variable importance methods have become so commonplace in predictive modeling that their role is often taken for granted. In standard linear regression problems, when interactions and higher order effects are absent, beta coefficients and $t$-statistics provide a concise and largely comprehensive summary of the influence of individual variables. However in more complicated black box models, the role of a predictor is dispersed and no single model parameter can quantify its importance. Leo Breiman introduced a groundbreaking approach for assessing feature importance in Random Forests, commonly known as the Permute-and-Predict (PaP) technique, which compares model accuracy with and without permuting a variable~\citep{Breiman}. This idea was later generalized into a model-agnostic framework known as model reliance, extending the same principle across a wide range of machine learning models~\citep{agnosticPerm}. As PaP gained widespread use, however, critics highlighted its potential limitations~\citep{Hooker, Gregorutti2017, Strobl}, prompting the development of alternative methods. One such approach is Leave-One-Covariate-Out (LOCO), an idea with roots in classical regression that was later formalized for machine learning~\citep{Lei, molnar2025}. LOCO involves dropping a variable, refitting, and comparing the resulting accuracy to that of the original model.

Although the idea of variable importance may seem straightforward, the term conceals fundamental ambiguity. Specifically, we must distinguish between the importance of a variable to a given predictive model and its importance to the prediction task itself~\citep{chenTrue}. Consider a feature that is strongly correlated with the response and thus receives a large regression coefficient, but is also highly collinear with another predictive feature. From a predictive-sensitivity perspective, the variable is important: perturbing it substantially alters predictions in a detrimental way. From a model-necessity perspective, however, the variable is dispensable: removing it causes little loss in performance because the correlated predictor compensates for its removal. This ambiguity underlies an important difference in classical regression metrics. Beta coefficients emphasize predictive sensitivity, whereas $t$-statistics emphasize model necessity. The sensitivity of the $t$-statistic to collinearity is well established and reflects its reduced statistical power in multiple regression \citep{HastieImp, DraperSmith1998,  Freedman1983}. The same distinction appears in machine learning metrics: PaP measures emphasize predictive sensitivity, while LOCO measures emphasize model necessity.

LOCO’s sensitivity arises because it retrains the model after removing or permuting the variable of interest. Thus, regardless of its strength in the initial model, if a variable’s predictive information is already captured by correlated predictors, then its removal has little effect on model errors, and LOCO deems it relatively unimportant. Through our derivations, we show that LOCO is approximately linearly (and inversely) related to the off-diagonal covariance structure of a variable transformation matrix. By contrast, PaP is unaffected by collinearity. Because the model is not retrained, permuting a variable inflates prediction error according to the magnitude of its coefficient, even when that variable is highly collinear with others. This distinction is not new. It is central to Chen et al.~\cite{chenTrue}, a fundamental component of Fisher et al.~\cite{agnosticPerm}, and closely related to the categorizations in Gregorutti et al.~\cite{Gregorutti2017} and Strobl et al.~\cite{Strobl} of marginal versus conditional importance. Within that nomenclature, our example feature is marginally important but not conditionally important, as its relationship to the response disappears when conditioned on collinear variables.

One might ask why Strobl et al.~\cite{strobl07} observed inflated importance for correlated predictors in Random Forests. The answer lies in the mechanics of tree-based methods, which differ from regression. Consider two highly correlated variables and two independent variables. The regression model could be
\[
\begin{aligned}
    y &= 2Cor_{1} + 2Cor_{2} + 2Ind_{3} + 2Ind_{4} + \epsilon\\ 
    &\approx 4Cor_{i}  + 2Ind_{3}  + 2Ind_{4} + \epsilon \text{ where } i \in \{1,2\}.
\end{aligned}
\]
In linear regression, this redundancy poses no problem as the effects of $Cor_1$ and $Cor_2$ are estimated jointly. In axis-aligned tree-based methods like Random Forests, however, splits are made greedily on a single variable at a time. Although random feature subsampling helps reduce this greediness, collinear predictors still gain elevated priority because their shared effect makes them strong candidates for early splits. This tendency to overstate importance also holds for strong negative correlations and coefficients. In contrast, when coefficients and correlations have opposite signs, tree-based methods understate the importance of correlated predictors \citep{BladenThesis}. 

In this paper we do not provide any arguments supporting the preference of one notion of variable importance over another. Instead, our primary motivation is to provide a clearer understanding of the precise differences in PaP and LOCO, particularly regarding their distinct sensitivities to collinearity. Specifically, we aim to provide a mathematical rationale for previously observed empirical findings by Bladen and Cutler: within a regression framework, PaP values exhibit a strong alignment with the regression coefficients, while LOCO values more closely resembled the $t$-statistics~\citep{BladenPerm}. To accomplish this, we introduce a highly general linear regression framework and show that, under certain assumptions, these empirical relationships necessarily follow from how each measure is calculated. We further support these derivations through a range of Monte Carlo simulations demonstrating their close alignment with empirical observation across a variety of conditions. Finally, although our derivations specifically quantify how PaP and LOCO operate within linear regression models, we conduct additional simulations to examine how well this theory generalizes to other models, such as Random Forests. Overall, by clarifying which data characteristics drive these variable importance measures, our results offer a clearer lens for interpreting variable importance in applied settings.

\section{Methods and Derivations} 
\label{s2:methods}

In this section, we provide derivations for the mathematical relationships between PaP, LOCO, $t$, and data attributes such as slope ($\boldsymbol{\beta}$), dimension ($p$), size ($n$), and collinearity (driven by $\Delta$).

\subsection{Latent Variable Framework}
\label{ss2:latent}

We begin by introducing a simple, but highly generalizable data structure that we can readily use in completing our derivations. Specifically, consider a traditional regression setup where the response variable $\mathbf{y}$ is expressed as an affine function of $\mathbf{X}$ and some error $\epsilon \overset{\mathrm{iid}}{\sim} \mathcal{N}(0, v = 0.1)$. To represent collinearity in the predictors, assume that $\mathbf{X}$ itself arises as a linear transformation of independent latent variables $\mathbf{Z} \sim \mathcal{N}(\mathbf{0}, \mathbf{I})$ via a correlation-inducing matrix $\mathbf{A}$. Substituting this structure into the regression model gives

\begin{equation}
\label{eq:reg_latent}
\underset{n\times 1}{\mathbf{y}}
= \underset{n\times p}{\mathbf{X}} \ \underset{p\times 1}{\boldsymbol{\beta}} + \underset{n\times 1}{\epsilon}
= \underset{n\times p}{\mathbf{Z}} \ \underset{p\times p}{\mathbf{A}} \ \underset{p\times 1}{\boldsymbol{\beta}} + \underset{n\times 1}{\epsilon},
\end{equation}

where

\[
\underset{p\times p}{\mathbf{A}} =\begin{bmatrix}
1 & \Delta & \cdots & \Delta \\
\Delta & 1 & \ddots & \vdots \\
\vdots & \ddots & \ddots & \Delta \\
\Delta & \cdots & \Delta & 1
\end{bmatrix}.
\]

%





For subsequent derivations, we can use the following formulation for $\mathbf{A}$:  $\mathbf{A}=\Delta\mathbf{J}+(1-\Delta)\mathbf{I}$, where $\mathbf{J}$ is a matrix of ones, $\mathbf{I}$ is the identity matrix, and $\Delta$ is a fixed real value in the range $[0, 1)$. This means
\begin{equation}
\label{eq:cov_x}
    cov(\mathbf{X})=\ cov(\mathbf{ZA})=\mathbf{A}cov(\mathbf{Z})\mathbf{A^T}=\mathbf{AA^T} =(2\Delta + (p-2)\Delta^2)\mathbf{J}+(1-\Delta)^2\mathbf{I}.
\end{equation}


\subsection{Derivation of Variable Importance for PaP}
\label{ss2:pap}

We begin by deriving the theoretical value for PaP.

\begin{definition}[$PaP$]
\label{def:pap}
The Permute-and-Predict variable importance measure can be defined as
\begin{equation}
\label{eq:pap_og}
    \text{PaP}_i = \sqrt{\mathscr{L}(\hat{f}(\mathbf{X}^v_{\pi(i)}),\mathbf{y}^v) - \mathscr{L}(\hat{f}(\mathbf{X}^v),\mathbf{y}^v)},
\end{equation}
where $\mathscr{L}$ is a selected loss function, $\hat{f}$ is a learned model, $\mathbf{X}^v$ is a validation predictor set, $\mathbf{y}^v$ is a validation response variable, and $\pi(i)$ denotes a permutation applied to the $i^{th}$ variable in the data set \cite{Breiman, agnosticPerm, Hooker, BladenPerm}.
\end{definition}

In plain terms, we are calculating the PaP importance for a given variable by permuting it in the validation set, generating new predictions, and comparing this new loss to the original loss. While many versions of PaP compute simple loss differences or compare root losses, we instead apply a square-root transformation to the raw difference, a monotonic adjustment motivated by observations in \cite{BladenPerm}.

Now consider the classical regression loss function $MSE$, which is minimized in linear regression and commonly used in the calculation of PaP. For simplicity, we assume $\mathbf{\hat{\boldsymbol{\beta}}} = \mathbf{\boldsymbol{\beta}}$. While this may seem a somewhat dubious assumption, any error in these regression coefficient estimates will inflate the $MSE$ in both the original and permuted models. Thus, as the variable importance measure is calculated as the difference in MSE terms, the affects of this assumption should be minimal - likely negligible - leaving the results both meaningful and robust for our purposes. 

From this assumption, it follows that the validation MSE of the original learned model simplifies to
\begin{equation}
\label{eq:imp_og}
    MSE(\hat{f}(\mathbf{X}^v), \mathbf{y}^v)=
    E\left[(\mathbf{\hat{y}}^v-\mathbf{y}^v)^2\right]= E\left[(\mathbf{X}^v\hat{\boldsymbol{\beta}}-\mathbf{X}^v\boldsymbol{\beta}+\mathbf{\mathbf{\epsilon}})^2\right]= E\left[\mathbf{\epsilon}^2\right]=\operatorname{Var}(\mathbf{\epsilon}).
\end{equation}

Next, suppose we implement a permutation to convert $\mathbf{X}^v$ to $\mathbf{X}^v_{\pi(i)}$. In this case, the learned coefficients $\mathbf{\hat{\boldsymbol{\beta}}}$ are unchanged. Because only the $i^{th}$ variable has been altered, $\hat{f}(\mathbf{X}^v_{\pi(i)})$ can be defined as
\[
\hat{f}(\mathbf{X}^v_{\pi(i)}) = \mathbf{X}^v\hat{\boldsymbol{\beta}}-\hat{\beta_i} \mathbf{x}^v_i + \hat{\beta_i} \mathbf{x}^v_{\pi(i)},
\]
where $\mathbf{x}^v_{\pi(i)}$ denotes a permutation of only the $i^{th}$ variable. Assuming $\beta=\hat{\beta}$, this yields
\begin{equation}
\label{eq:pap_perm}
    \begin{split}
        MSE(\hat{f}(\mathbf{X}^v_{\pi(i)}), \mathbf{y}^v) &=
        E\left[(\hat{f}(\mathbf{x}^v_{\pi(i)})-\mathbf{y}^v)^2\right]\\
        &= E\left[\left( (\mathbf{X}^v\hat{\boldsymbol{\beta}}-\hat{\beta_i} \mathbf{x}^v_i + \hat{\beta_i} \mathbf{x}^v_{\pi(i)}) - \mathbf{X}^v\boldsymbol{\beta}-\mathbf{\epsilon}\right)^2\right] \\
        &= E\left[\left( -\hat{\beta_i} \mathbf{x}^v_i + \hat{\beta_i} \mathbf{x}^v_{\pi(i)} - \mathbf{\epsilon}\right)^2\right] \\
        &= E[\mathbf{\epsilon^T\epsilon}]+ \hat{\beta_i}^2\left( E[{\mathbf{x}^v_i}^T \mathbf{x}^v_i]+E[{\mathbf{x}^v_{\pi(i)}}^T \mathbf{x}^v_{\pi(i)}] \right)\\
        &= \operatorname{Var}(\epsilon)+2\hat{\beta_i}^2 \operatorname{Var}(\mathbf{x}^v_i),
    \end{split}
\end{equation}
since $\mathbf{\epsilon}$, $\mathbf{x}^v_i$, and $\mathbf{x}^v_{\pi(i)}$ are all zero-mean and independent of each other, allowing their means and cross terms to be removed. Substituting \cref{eq:imp_og} and \cref{eq:pap_perm} into \cref{eq:pap_og} yields the closed-form solution in \cref{th:pap}.

\begin{theorem}\label{th:pap}
Under the latent variable framework in \cref{eq:reg_latent} and assuming $\beta=\hat{\beta}$, the PaP importance has the closed-form
\begin{equation}
\label{eq:pap}
    \text{PaP}_i= \beta_i\sqrt{2\operatorname{Var}(\mathbf{x}^v_i)} = \beta_i\sqrt{2(1+ (p-1)\Delta^2)}
\end{equation}
where $\operatorname{Var}(\mathbf{x}^v_i)$ is found by substituting a value of $1$ into $\mathbf{J}$ and $\mathbf{I}$ in \cref{eq:cov_x}. \hfill $\textcolor{qedblue}\blacksquare$
\end{theorem}

Thus, in the proposed regression framework, $\text{PaP}_i$ depends only on the magnitude of the coefficient $\beta_i$ and the variance of the feature. Importantly, the degree of collinearity (defined by $\Delta$) has no direct effect on $\text{PaP}_i$ aside from its influence on $\operatorname{Var}(\mathbf{x}^v_i)$. Consequently, $\text{PaP}_i$ is independent of collinearity as long as there are no co-occurring effects on variable scale. This independence can also be guaranteed when variables are standardized (S.T. $\operatorname{Var}(\mathbf{x}^v_i) = 1,  \text{for all } i$) \cite{Hooker}. We also note that a Drop-and-Predict technique admits an essentially identical theoretical formulation to PaP, differing only by a multiplicative factor of 2 \cite{BladenPerm}.

\subsection{Derivation of Variable Importance for LOCO}
\label{ss2:loco}

Similar to the expression for PaP, the leave-one-covariate-out variable importance measure can be expressed in terms of loss functions.

\begin{definition}[$LOCO$]
\label{def:loco}
The leave-one-covariate-out variable importance measure can be defined as
\begin{equation}
\label{eq:loco_og}
    \text{LOCO}_i = \sqrt{\mathscr{L}(\hat{f}^t_{\mathbf{\pi(i)}}(\mathbf{X}^v),\mathbf{y}^v) - \mathscr{L}(\hat{f}(\mathbf{X}^v),\mathbf{y}^v)},
\end{equation}
where the new notation ($\hat{f}^t_{\mathbf{\pi(i)}}$) represents a model retrained on the training data with the $i^{th}$ variable permuted: $\mathbf{X}^t_{\pi(i)}$ \cite{Lei, Hooker, BladenPerm}.
\end{definition}

In other words, we calculate the LOCO importance for a given variable by permuting it in the training set, refitting the model, generating new validation predictions, and comparing this new loss to the original loss. We again take the square root to place this on a scale familiar to regression metrics \citep{BladenPerm}.

When a training variable is permuted and a model is refit, the result is a different set of coefficients $\hat{f}^t_{\mathbf{\pi(i)}}= \hat{\boldsymbol{\beta}}_{\mathbf{\pi(i)}}$. For simplicity we assume that the value of the permuted coefficient is equally absorbed into the coefficients of all other attributes which it is correlated with. We express this as
\begin{equation}
\label{eq:newb}
\hat{\boldsymbol{\beta}}_{\mathbf{\pi(i)}} = 
\begin{bmatrix}
\hat{\beta}_1' \\
\vdots \\
\hat{\beta}_{i-1}' \\
\hat{\beta}_i'=0 \\
\hat{\beta}_{i+1}'\\
\vdots \\
\hat{\beta}_{p}'
\end{bmatrix}=
\begin{bmatrix}
\hat{\beta}_1 +\hat{\beta}_i*c   \\
\vdots \\
\hat{\beta}_{i-1} +\hat{\beta}_i*c  \\
\hat{\beta}_i - \hat{\beta}_i \\
\hat{\beta}_{i+1}+\hat{\beta}_i*c \\
\vdots \\
\hat{\beta}_{p}+\hat{\beta}_i*c 
\end{bmatrix} = 
\hat{\boldsymbol{\beta}}+\boldsymbol{\delta},
\end{equation}

where $c$ is a coefficient controlling how much of the original magnitude is absorbed (hereafter called the coefficient of absorption) and  $\boldsymbol{\delta}^i = \begin{bmatrix} \hat{\beta}_i*c, \dots, \hat{\beta}_i*c, -\hat{\beta}_i, \hat{\beta}_i*c, \dots, \hat{\beta}_i*c \end{bmatrix} ^T$. 
Note that the equal absorption of the permuted coefficient is supported by the empirical results in \cref{ss3:c}. Also note that the coefficient of absorption, which we will discuss in greater detail later on, is heavily dependent on the degree of collinearity in the data.

Now we use the results of \cref{eq:newb} to derive the new loss function:
\begin{equation}
\label{eq:loco_perm}
    \begin{split}
        MSE(\hat{f}^t_{\pi(i)}(\mathbf{X}^v), \mathbf{y}^v)&=
        E\left[(\hat{f}^t_{\pi(i)}(\mathbf{X}^v)-\mathbf{y}^v)^2\right]=E\left[\left(\mathbf{X}^v(\hat{\boldsymbol{\beta}}+\boldsymbol{\delta}) - \mathbf{X}^v\boldsymbol{\beta}-\mathbf{\epsilon}\right)^2\right] \\
        &=E\left[\left(\mathbf{X}^v \boldsymbol{\delta} - \mathbf{\epsilon}\right)^2\right] = E\left[\left(\mathbf{Z}^v \mathbf{A}\boldsymbol{\delta} - \mathbf{\epsilon}\right)^2\right] \\
        &=E[\mathbf{\epsilon^T\epsilon}]+ E\left[\boldsymbol{\delta}^T\mathbf{A}^T{\mathbf{Z}^v}^T \mathbf{Z}^v \mathbf{A}\boldsymbol{\delta}\right] = \operatorname{Var}(\epsilon)+ E\left[\boldsymbol{\delta}^T\mathbf{A}^T\mathbf{A}\boldsymbol{\delta}\right].
        \end{split}
        \end{equation}

Substituting \cref{eq:cov_x} for $\mathbf{A}^T\mathbf{A}$ in \cref{eq:loco_perm} gives
\begin{equation}
\label{eq:loco_perm2}
    \begin{split}
        MSE(\hat{f}^t_{\pi(i)}(\mathbf{X}^v), \mathbf{y}^v) &= \operatorname{Var}(\epsilon)+ E\left[\boldsymbol{\delta}^T((2\Delta + (p-2)\Delta^2)\mathbf{J}+(1-\Delta)^2\mathbf{I})\boldsymbol{\delta}\right]\\
        &= \operatorname{Var}(\epsilon)+ \hat{\beta}_i^2[(2\Delta + (p-2)\Delta^2)((p-1)c - 1)^2 + (1-\Delta)^2(1+(p-1)c^2)].
    \end{split}
\end{equation}

Since $\mathbf{\epsilon}$ and $\mathbf{X}^v$ are zero-mean and independent of each other, their means and cross terms vanish. Substituting the results of \cref{eq:imp_og} and \cref{eq:loco_perm2} into \cref{eq:loco_og} yields a closed-form solution for LOCO (stated in \cref{th:loco_exact}).

\begin{theorem}\label{th:loco_exact}
Under the framework in \cref{eq:reg_latent} and assuming $\beta=\hat{\beta}$,
the LOCO importance has the closed-form
\begin{equation}
\label{eq:loco_exact}
    \text{LOCO}_i=\beta_i\sqrt{(1-\Delta)^2(1+(p-1)c^2) + (2\Delta + (p-2)\Delta^2)((p-1)c - 1)^2}.
    \rlap{\hspace{0.92in}$\textcolor{qedblue}\blacksquare$}
\end{equation}
\end{theorem}

While $c$ was initially unknown, being learned from model metrics rather than a fixed parameter, it follows that theoretical $c$ should be a function of $\Delta$ and $p$. 
\begin{conjecture}[$c$]
\label{conj:c}
Based on empirical observations, we postulate that $c$ takes the following form:
\begin{equation}
\label{eq:c}
    c(\Delta, p) = \frac{2\Delta + (p-2)\Delta^2}{1-4\Delta + 2p\Delta  + p^2\Delta^2 -3p\Delta^2 + 3\Delta^2}.
\end{equation}
\end{conjecture}
Note that \cref{conj:c} is formulated by assuming the relationship subsequently shown in \cref{th:loco_t} and then deriving the exact $c$ value that yields this relationship. However, we will show that (1) empirical $c$ values fit this conjecture very closely (see \cref{ss3:c}), and (2) if this conjecture is true, then \cref{th:loco_simp} and \cref{th:loco_t} necessarily follow.





\cref{eq:c} can be easily evaluated at the boundary points of $\Delta$. If $\Delta=0$, then all variables in $\mathbf{X}$ are independent and $c = 0$, meaning there is no absorption of the permuted coefficient. In this case, substituting $\Delta = 0$ and $c = 0$ into \cref{eq:loco_exact} yields $\text{LOCO}_i = \hat{\beta_i}$. If $\Delta=1$, then all variables in $\mathbf{X}$ are perfectly collinear and $c =\frac{1}{p-1}$. In this situation, substituting these into \cref{eq:loco_exact} yields $\text{LOCO}_i= 0$.

By substituting \cref{eq:c} into \cref{eq:loco_exact} and simplifying (see \cref{app:loco_c}), we obtain the following representation of LOCO.

\begin{theorem}\label{th:loco_simp}
Given \cref{eq:loco_exact}, LOCO simplifies via \cref{conj:c} to
\begin{equation}
\label{eq:loco_c_simp}
\text{LOCO}_i = \beta_i (1-\Delta)\sqrt{\frac{(1 + (p-1)\Delta)^2}{(1 + (p-2)\Delta)^2 + (p-1)\Delta^2}} = \beta_i (1-\Delta)\sqrt{1 + c}.
\rlap{\hspace{1in}$\textcolor{qedblue}\blacksquare$}
\end{equation}
\end{theorem}

A very natural approximation for this importance is 
\begin{equation}
\label{eq:loco_approx}
    \text{LOCO}_i \approx \beta_i(1-\Delta),
\end{equation}
since $1 < \sqrt{1+c} < \frac{p}{p-1}$. Put simply, these formulations show that collinearity systematically shrinks LOCO, in contrast to PaP, which depends only on predictor variability.

While the theoretical value of LOCO is well defined in \cref{eq:loco_c_simp}, Bladen and Cutler showed that it has a close connection to regression $t$-statistics~\cite{BladenPerm}. 

\begin{definition}[$t$-statistic]
\label{def:t}
The $t$-statistic for coefficient $\hat{\beta}_i$ is defined as
\begin{equation}
\label{eq:t_stat_og}
    t_i = \frac{\hat{\beta}_i}{SE(\hat{\beta}_i)} = \frac{\hat{\beta}_i}{\sqrt{\operatorname{Var}(\epsilon)(\mathbf{X^{T}X})^{-1}_{ii}}}.
\end{equation}
\end{definition}

\begin{theorem}[$t$-statistic as function of $\Delta$]\label{th:t_stat}
Under the framework in \cref{eq:reg_latent} and assuming $\beta=\hat{\beta}$, the $t$-statistic of \cref{eq:t_stat_og} can be expressed as
\begin{equation}
\label{eq:t_stat}
    t_i = \beta_i(1-\Delta) \sqrt{\frac{(n-1)(1+(p-1)\Delta)^2}{\operatorname{Var}(\epsilon)[(1 + (p-2)\Delta)^2 + (p-1)\Delta^2]}}.
    \rlap{\hspace{1.59in}$\textcolor{qedblue}\blacksquare$}
\end{equation}
\end{theorem}

A complete mathematical proof deriving \cref{th:t_stat} from \cref{eq:t_stat_og} is available in \cref{app:t_stat}.

\begin{theorem}\label{th:loco_t}
Given \cref{eq:loco_c_simp}, \cref{eq:t_stat}, and \cref{conj:c}, the precise relationship between the $t$-statistic and LOCO can be mathematically quantified as
\begin{equation}
\label{eq:loco_t}
    t_i = \text{LOCO}_i\sqrt{\frac{n-1}{\operatorname{Var}(\epsilon)}}.
    \rlap{\hspace{2.64in}$\textcolor{qedblue}\blacksquare$}
\end{equation}
\end{theorem}

\subsection{Simulation Design}
\label{ss2:sim_design}

We now offer a system of simulations to show the correctness and relevance of our resulting derivations and to validate our assumptions. The general setup aligns with the latent variable framework introduced in \cref{ss2:latent}.

\subsubsection{Data Generation}
\label{ss2:data}

We simulate $\mathbf{X}$ as a set of $p$ predictor variables with $n$ observations. These are constructed by applying the linear transformation $\mathbf{A} = \Delta\mathbf{J} + (1 - \Delta)\mathbf{I}$ to a matrix $\mathbf{Z}$ of independent standard normal variables, where $\mathbf{J}$ is the $p \times p$ matrix of ones and $\Delta \in [0, 0.99)$ controls the correlation structure. The response variable $\mathbf{y}$ is defined as a linear function of $\mathbf{X}$ with identical coefficients $\boldsymbol{\beta}$ and additive noise $\epsilon \overset{\mathrm{iid}}{\sim} \mathcal{N}(\operatorname{Mean} = 0, \operatorname{Var} = 0.1)$.

We vary $\Delta \in [{0, 0.11, 0.22, \dots, 0.88, 0.99}]$, $p \in [{3, 6, 9, 12}]$, and $n \in [{20, 63, 200, 632, 2000}]$ systematically. For each parameter combination, we generate 100 Monte Carlo iterations using identically structured training and validation sets.

\subsubsection{Linear Model Implementation}

For each simulation, we first fit a linear regression model to the training set and evaluate the metrics $\text{PaP}_i$, $\text{LOCO}_i$, $t_i$, and $c$ on the validation set. These empirically observed quantities are then compared to their corresponding theoretical or assumed expressions in \cref{eq:pap}, \cref{eq:loco_c_simp}, and \cref{eq:c} to assess the validity of our derivations under the assumed linear data-generating process. Comparisons are made by computing the relative difference: $\frac{\text{empirical} - \text{theory}}{\text{theory}}$. Functional relationships are also assessed by plotting both the empirical and theoretical trends across $\Delta$.

\subsubsection{Random Forest Implementation}

To assess the robustness of our results under model misspecification and nonlinearity, we repeat the full set of simulations using a Random Forest model in place of the linear predictor. The data generation and evaluation settings remain identical. For each simulation, we fit a Random Forest on the training set and compute the same set of metrics on the validation set. We fix the number of candidate variables at each split to $mtry = \max\left(2,~\text{floor}(\frac{p}{3})\right)$ to maintain consistency across values of $p$. This nonparametric modeling approach allows us to test whether the theoretical insights derived under linear assumptions retain utility in more flexible model classes.

Comparisons are again made by computing the relative difference for PaP ($\frac{\text{empirical} - \text{theory}}{\text{theory}}$) or the absolute difference for LOCO ($\text{empirical} - \text{theory}$) to accommodate the noise associated with Random Forest importances and difficulty of performing a relative comparison for small theoretical values. Functional relationships are also assessed by plotting the empirical and theoretical trends across $\Delta$.

\section{Results} 
\label{s:results}

We employ \texttt{R} \citep{R-base} and the \texttt{tidyverse} \citep{tidyverse} to carry out the simulations defined in \cref{ss2:sim_design}. Random Forest analyses are implemented using the \texttt{randomForest} package \citep{randomForest}. For visualization, we extract a meaningful subset of our results and plot comparisons using \texttt{ggplot2} \citep{ggplot2} together with \texttt{ggh4x} \citep{ggh4x}. 

In \cref{ss3:lin_imp}, we provide plots comparing the means and spreads of empirical linear regression metrics with the theoretical values derived in \cref{eq:pap} and \cref{eq:loco_c_simp} (see \cref{app:t_results} for empirically testing \cref{eq:loco_t}). We also offer visualizations of the functional relationship between these metrics and $\Delta$. These visualizations are also applied to show how empirical Random Forest metrics compare to the theoretical values and to validate our assumptions for $c$ in \cref{eq:c}.

\subsection{Linear Model Importance Results}
\label{ss3:lin_imp}
\begin{figure}[htbp]
\begin{center}
\centerline{\includegraphics[width=0.9\columnwidth]{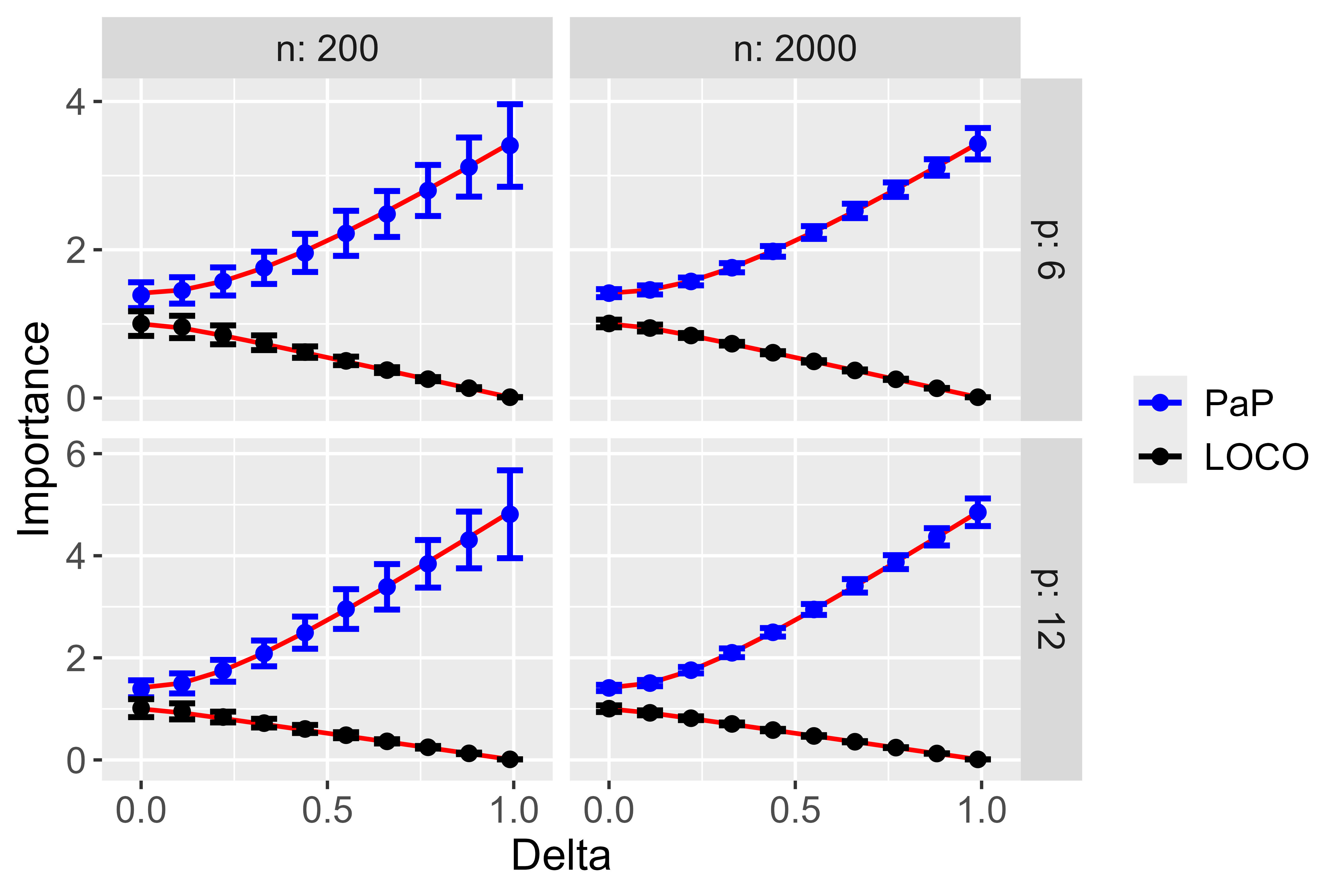}}
\caption{Mean and 2 Standard deviation error bars for assessing the relationship between importance values and $\Delta$ (see \cref{eq:pap} and \cref{eq:loco_c_simp}) for data defined in \cref{ss2:data}. Red lines denote theoretical projections. Plots show importance values for raw data and are faceted by $p$ and $n$ dimensions. Plots show that PaP increases with $\Delta$, while LOCO has a nearly linear negative relationship with $\Delta$ (See \cref{eq:loco_approx}).}
\label{fig:cov_delta}
\end{center}
\vskip -0.2in
\end{figure}

We begin by conducting parity comparisons for importance in linear models. The initial question tested by the simulation is how well the derived theoretical expressions match the empirical results for both LOCO and PaP as a function of $\Delta$. \cref{fig:cov_delta} illustrates the average empirical values computed for each measure with 2SE error bars as a function of $\Delta$ for sample values of $n$ and $p$. Recall from \cref{eq:pap} that $\text{PaP}_i= \beta_i\sqrt{2\operatorname{Var}(\mathbf{x}^v_i)}$, and for our simulations, $\operatorname{Var}(\mathbf{x}^v_i) = 1+(p-1)\Delta^2$. As a result, both the theoretical and empirical values of PaP increase with $\Delta$ under this covariance structure found in \cref{fig:cov_delta}. In contrast, LOCO decreases nearly linearly with increasing $\Delta$, closely matching the relationship identified by \cref{eq:loco_c_simp}.

\subsection{Random Forest Importance Results}
\begin{figure}
\begin{center}
\centerline{\includegraphics[width=0.9\columnwidth]{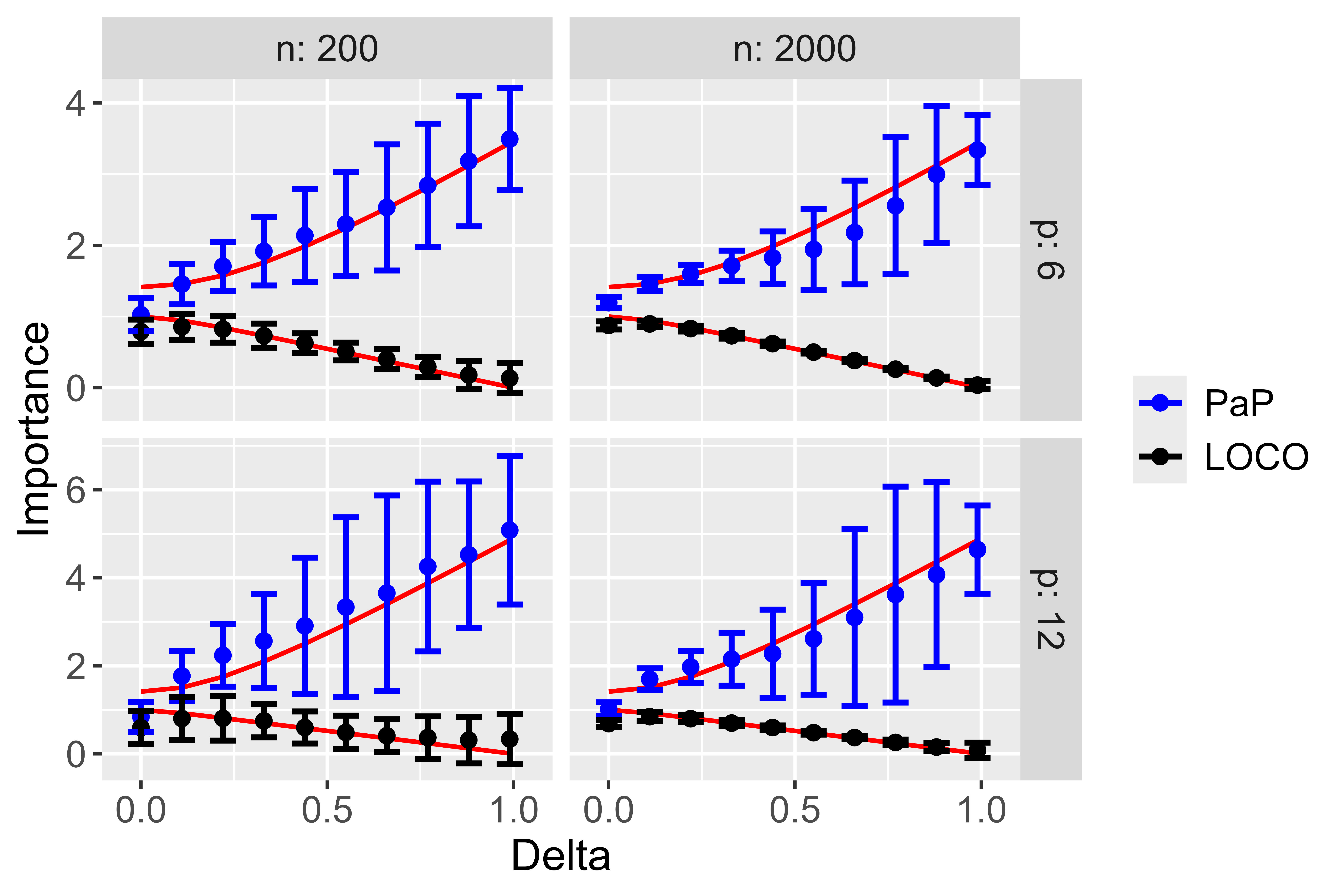}}
\caption{Mean and 2 Standard deviation error bars for assessing the relationship between Random Forest importance values and $\Delta$ for data defined in \cref{ss2:data}. Red lines denote theoretical projections for linear data. Plots show importance values for raw data and are faceted by $p$ and $n$ dimensions. Plot trends appear similar to \cref{fig:cov_delta}. However, there are clear anomalies regarding bias within the plot, especially for $\Delta = 0$.}
\label{fig:rf_cov_delta}
\end{center}
\vskip -0.2in
\end{figure}

\cref{fig:rf_cov_delta} illustrates how the distribution of empirical Random Forest importance values, the theoretical expectations, and the relationship between them behave as functions of $\Delta$. Theoretical values are again derived from \cref{eq:pap} and \cref{eq:loco_c_simp}. When compared to the near-perfect agreement between empirical and theoretical values for the linear model in \cref{fig:cov_delta} (aside from some small-sample error), the Random Forest results in \cref{fig:rf_cov_delta} are notably less consistent.


Empirical patterns for both PaP and LOCO generally align with the theoretical curves, though some systematic deviations emerge. For PaP, most differences are reasonably modest, except at $\Delta = 0$, where empirical values fall notably below the theoretical expectation. A similar anomaly is observed for LOCO at this point. Furthermore, LOCO exhibits a more consistent bias pattern: empirical values tend to lie below the theory for small $\Delta$ and above it for larger $\Delta$.  Overall, despite some of these deviations, the trends show that both PaP and LOCO align closely with our theoretical predictions, particularly for larger $n$. This is indicative that the collinearity effects quantified in our framework extend and provide meaningful insight into variable importance for Random Forests and likely other complex machine learning regression models.

\subsubsection{\texorpdfstring{$c$}{c} Results}
\label{ss3:c}

\begin{figure}[htbp]
\begin{center}
\centerline{\includegraphics[width=0.9\columnwidth]{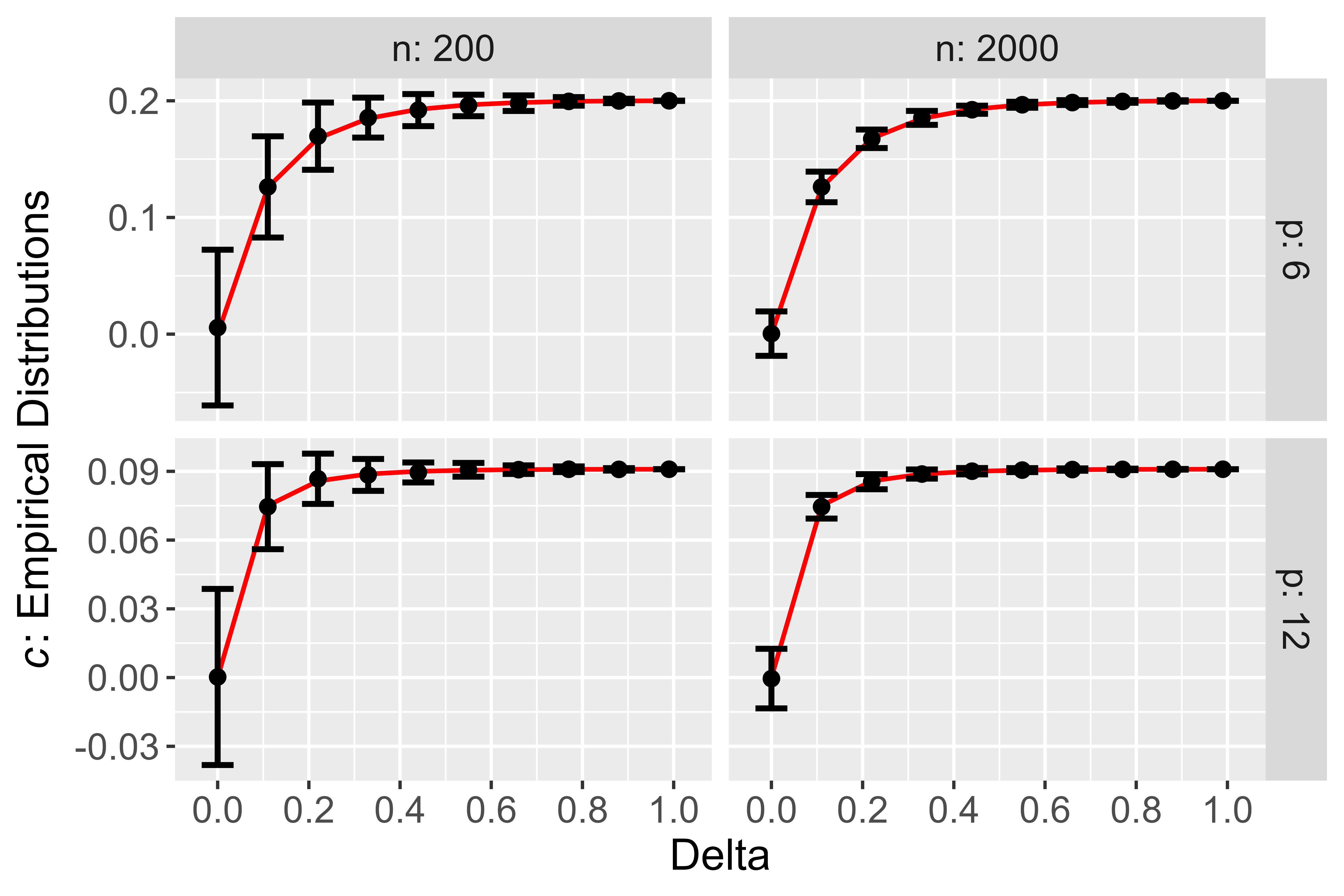}}
\caption{Mean and 2 Standard deviation error bars for assessing the relationship between c values and $\Delta$ (see \cref{conj:c}) for data defined in \cref{ss2:data}. 
Red lines denote theoretical projections. Plots are faceted by $n$ and $p$ dimensions. Plots show increased stability as $n$ increases, $\Delta$ increases, and $p$ decreases. Plots also show sublinear growth converging toward values of $\frac{1}{p-1}$. Plots are free of any clear bias.}
\label{fig:c_delta}
\end{center}
\vskip -0.2in
\end{figure}

\cref{fig:c_delta} illustrates the relationship between $\Delta$, $p$, and the functional value of $c$ in a linear regression model. It also visualizes the relative differences between the empirical distribution and the proposed value of $c$, postulated in \cref{conj:c}. As shown, $c$ increases monotonically with $\Delta$, though the increase is clearly sublinear and converges toward an upper bound of $\frac{1}{p-1}$. This limit is highly intuitive: when one variable is removed or permuted, the remaining $p-1$ variables absorb the effect, and if that influence is distributed uniformly, each receives a share of $\frac{1}{p-1}$. \cref{fig:c_delta} shows an intuitive decrease in variance as $n$ increases and a more surprising decrease in variance with increasing $p$ and $\Delta$. Importantly, no observable bias is present in the plot, visually suggesting that \cref{eq:c} provides an unbiased estimate of $c$ and is therefore highly relevant for the LOCO derivations presented in \cref{ss2:loco}.

\subsubsection{Small Sample Bias}
\label{s3:ss_bias}

\begin{figure}[htbp]
\begin{center}
\centerline{\includegraphics[width=0.9\columnwidth]{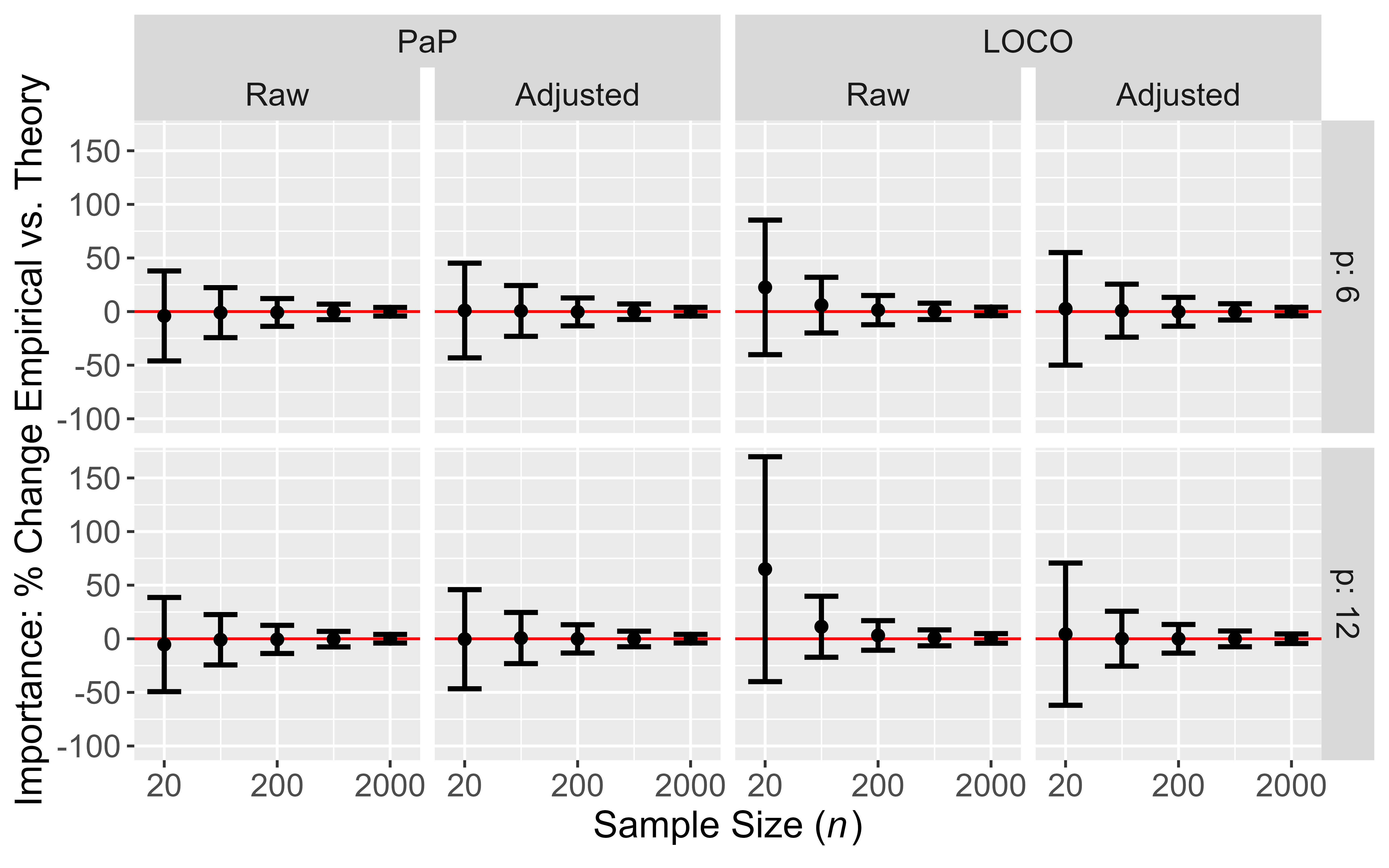}}
\caption{Mean and 2 Standard deviation error bars for comparing empirical importance values to the derived theoretical importances (see \cref{eq:pap} and \cref{eq:loco_c_simp}) for data defined in \cref{ss2:data}. Plots show sample size ($n$) on the x-axis and are faceted by dimensions ($p$), importance metric, and whether a sample sample adjustment was applied to the theoretical values. Plots show increased stability as $n$ increases and as $p$ decreases. Plots also show a small sample bias for both raw metrics, but particularly for LOCO. This bias is especially noticeable for small $n$ and large $p$ combinations. Conversely, the adjusted metrics are free of any clear visual bias.}
\label{fig:combo_parity}
\end{center}
\vskip -0.2in
\end{figure}

\cref{fig:combo_parity} illustrates the relative differences between empirical and theoretical importance values for the linear regression model. As expected, the variance decreases with increasing $n$ and decreasing $p$. A subtle small-sample bias can be observed for the raw PaP metrics, where theoretical values slightly exceed empirical estimates. Meanwhile, raw LOCO values exhibit a clear upward bias, with empirical values exceeding theory—most notably when $n$ is small and $p$ is large.

To address these small-sample biases, we apply a sample size correction to our evaluation metrics. For the PaP metric, we adopt Bessel’s correction — a widely used adjustment in statistical estimation — most notably in the unbiased estimation of sample variance \citep{kenney_bias, stigler_bias, reichel_bias}. This correction typically takes the form $\frac{n}{n-1}$, or its inverse $\frac{n-1}{n}$, depending on the formulation. Incorporating this adjustment into \cref{eq:pap} yields:
\begin{equation}
\label{eq:pap_bias}
    \text{PaP}_i= \beta_i\sqrt{2\operatorname{Var}(\mathbf{x}^v_i)}\times\frac{n-1}{n}.
\end{equation}

LOCO is addressed with a different, yet similarly well-established, correction drawn from regression modeling. This finite-sample bias correction factor, often referred to as a degrees-of-freedom adjustment, takes the form $\sqrt{\frac{n}{n-p}}$ or its inverse $\sqrt{\frac{n-p}{n}}$, depending on the desired scale and interpretation \citep{efron_bias, shao_bias}. This adjustment especially compensates for overfitting and bias that can occur when the number of predictors $p$ approaches the sample size $n$. Incorporating this term into \cref{eq:loco_c_simp} yields:
\begin{equation}
\label{eq:loco_bias}
        \text{LOCO}_i= \beta_i(1-\Delta)\sqrt{1 + c}\times\sqrt{\frac{n}{n-p}}.
\end{equation}

\cref{fig:combo_parity} shows that implementing these adjustments retains the general variance patterns of the raw metrics. However, small-sample biases are effectively eliminated, so that empirical distributions are centered closely around their theoretical counterparts.

\section{Discussion and Conclusion}

Our results provide new theoretical insights into two widely used variable importance measures: Permute-and-Predict and Leave-One-Covariate-Out. By deriving formal mathematical expressions for each, we show how they are precisely linked to closed-form regression metrics and data characteristics of dimension ($p$), collinearity ($\Delta$), and effect size ($\beta_i$). This provides a more grounded understanding of what these measures reveal about the structure of the data itself, especially where collinearity is concerned. Empirical simulations show that the proposed framework provides a very precise model for how these variable importance measures react to collinearity in linear regression and generalizes very well to other models such as Random Forests. Placed in the context of prior work, our results confirm earlier assessments of importance methods by \cite{Hooker, chenTrue}, extend prior comparisons to regression metrics by \cite{BladenPerm}, and solidify the theoretical foundations for interpreting variable importance in relation to data characteristics.

The conceptual interpretation of the derived relationships can be summarized as follows. If a variable captures information about the response that could also be represented by other predictors, the PaP measure assigns it importance regardless of how little unique information it contributes. In contrast, LOCO assigns little importance to such variables, since retraining allows that information to be absorbed by other predictors. In the framework of this paper, the degree of shared information is captured by collinearity, and both PaP and LOCO can be expressed precisely in terms of the true coefficient $\beta_i$ and, for LOCO, the degree of collinearity ($\Delta$). In more general settings where variables may share nonlinear mutual information, these derivations are not directly applicable. Nevertheless, nonlinear relationships among predictors would likely produce similar effects on these measures, provided that the regression model has sufficient capacity to map the different representations of information to the response.

\subsection{Limitations and Considerations}

Despite these contributions, several limitations should be acknowledged. Our derivations are framed in standard regression settings and, while they provide valuable and legitimate insights, some extensions to more complex or nonlinear data structures remain to be formally developed. Note that while we tested these derivations in the context of Random Forests, the experiments still assumed the underlying linear relationships and collinearity structure defined in \cref{ss2:latent}. Moreover, while Random Forests are an illustrative test case to explore these derivations for non-parametric models, further evaluation in additional contexts is warranted before generalizing broader conclusions.

Additionally, the finite-sample corrections considered are grounded in classical regression theory, so their applicability to high-dimensional data or strongly regularized models requires further investigation. Our empirical analysis also has some limitations as it is necessarily constrained by the chosen simulation designs, which reflect particular distributions and dimensionalities. While no single study can capture every possible condition, these results establish a solid foundation that future work can build on to explore broader settings and applications. We are optimistic that these contributions will not only clarify existing methods but also inspire new directions in the development of interpretable machine learning tools and theories.

\subsection{Extensions and Applications}

Although our analysis was carried out in an idealized regression framework, the results naturally extend to other common data structures. One clear example involves $K$ independent blocks of correlated variables. When applying our theory to such cases, the effective dimensionality ($p_k$) corresponds to the correlated block which a variable belongs to. We note that prior research uses a dataset where one such block contains variables that are independent of each other as well as the rest of the blocks \citep{Strobl}. In this case, our theoretical results apply to each block, with the independent block characterized by $\Delta_{i,j}=0$ where $i$ or $j$ is in the indices of block $K$. These extensions follow directly from the logic of our derivations and illustrate the generality of the framework. Time-series data structures represent a more challenging setting. While the exact results presented here would not hold, the general trends remain informative: both PaP and LOCO are shaped by the degree of collinearity in the data, and so the insights provided here should still serve as a reasonable approximation. Exploring this connection more formally offers an interesting avenue for future work.

Additionally, extending these derivations to alternative importance frameworks (e.g., Global Shap values \citep{shap}, Partial Dependence Plot importances \citep{pdpvip}) could provide further clarity and help unify interpretability methods across modeling paradigms \cite{verdinelli2024feature}. More extensive simulation studies and applications to real-world datasets would further strengthen the validation of the theoretical derivations of PaP and LOCO, particularly in modern machine learning techniques. Beyond theoretical contributions, these results offer practical guidance by enabling importance measures to be understood and interpreted with greater confidence. Taken together, this work provides theoretical derivation bridged with empirical validation, laying a strong foundation for more principled and reliable use of variable importance in applied machine learning.

\section*{Code Access}
The code used to generate the data and figures in this article can be found at:
\url{https://github.com/KelvynBladen/theoryVarImportance}

\appendix

\section{Derivation of LOCO Variable Importance}
\label{app:loco_c}

We now offer a mathematical proof for the simplified version of LOCO shown in \cref{th:loco_simp}.

\begin{proof}
To prove \cref{th:loco_simp} by substituting the assumed value of $c$ from \cref{eq:c} into \cref{eq:loco_exact}, we begin by recalling and restating their mathematical expressions:
\begin{equation}
\tag{\ref{eq:loco_exact}}
    \text{LOCO}_i=\beta_i \sqrt{(1-\Delta)^2(1+(p-1)c^2)+ (2\Delta + (p-2)\Delta^2)((p-1)c - 1)^2}
\end{equation}

and
\begin{equation}
\tag{\ref{eq:c}}
    c = \frac{2\Delta + (p-2)\Delta^2}{1-4\Delta + 2p\Delta  + p^2\Delta^2 -3p\Delta^2 + 3\Delta^2}.
\end{equation}

For simplicity, we define the denominator as $d = 1-4\Delta + 2p\Delta  + p^2\Delta^2 -3p\Delta^2 + 3\Delta^2$. We then split the summand of \cref{eq:loco_exact} into terms: $s_1 = (1-\Delta)^2(1+(p-1)c^2)$ and $s_2 = (2\Delta + (p-2)\Delta^2)((p-1)c - 1)^2$. We will further split these into factors containing $c$: $s_1' = (p-1)c^2 + 1$ and $s_2' = (p-1)c - 1$. 

We'll simplify both of the factors, starting with $s_2'$:
\begin{equation}
\label{eq:s2p}
\begin{split}
    s_2' &= (p-1)c - 1 = \frac{(p-1)(2\Delta + (p-2)\Delta^2) - d}{d}\\
    &= \frac{2p\Delta + p^2\Delta^2 -2p\Delta^2 - 2\Delta - p\Delta^2+2\Delta^2-1+4\Delta - 2p\Delta - p^2\Delta^2 + 3p\Delta^2 - 3\Delta^2}{d}\\
    &= \frac{-1 + 4\Delta - 2\Delta +2\Delta^2 - 3\Delta^2 + 2p\Delta - 2p\Delta + 3p\Delta^2 - 2p\Delta^2 - p\Delta^2 + p^2\Delta^2 - p^2\Delta^2}{d}\\
    &=\frac{-1 + 2\Delta  - \Delta^2}{d}= \frac{- (1-\Delta)^2}{d}.
\end{split}
\end{equation}

Squaring $s_2'$ and multiplying it by $2\Delta + (p-2)\Delta^2$ yields the second term:
\begin{equation}
\label{eq:s2}
    s_2 = \frac{(1-\Delta)^4(2\Delta + (p-2)\Delta^2)}{d^2}= \frac{(1-\Delta)^4c}{d}.
\end{equation}

Now we alter $s_1'$, using the relationship found in \cref{eq:s2p}:
\begin{equation}
\label{eq:s1p}
\begin{split}
    &s_1' = (p-1)c^2 + 1\\
    &= (p-1)c^2 - c + c +1\\
    &= ((p-1)c - 1)c + 1 + c\\
    &=\frac{- (1-\Delta)^2 c}{d} + 1 + c.
\end{split}
\end{equation}

Multiplying $s_1'$ by $(1-\Delta)^2$ yields:
\begin{equation}
\label{eq:s1}
    s_1 = \frac{- (1-\Delta)^4 c}{d} + (1-\Delta)^2 (1 + c).
\end{equation}

Now, we have each of the terms for the summand in \cref{eq:loco_exact}. Summing \cref{eq:s1} and \cref{eq:s2} and simplifying yields:
\begin{equation}
\label{eq:summand}
     s_1 + s_2 = \frac{(1-\Delta)^4 c}{d} + \frac{- (1-\Delta)^4 c}{d} + (1-\Delta)^2(1 + c) = (1-\Delta)^2(1 + c).
\end{equation}

Finally, inserting \cref{eq:summand} back into \cref{eq:loco_exact} yields:
\begin{equation}
\tag{\ref{eq:loco_c_simp}}
    \mathbf{LOCO_i} = \beta_i\sqrt{s_1 +s_2} = \beta_i\sqrt{(1-\Delta)^2(1 + c) }=\beta_i(1-\Delta)\sqrt{1 + c}. \qed
\end{equation}
\end{proof}

\section{Derivation of \texorpdfstring{$t$}{t}-statistic Formula}
\label{app:t_stat}

We now offer a mathematical proof for \cref{th:t_stat} where the $t$-statistic is expressed as a function of $\hat{\beta}_i$, $p$, and $\Delta$. Proving \cref{th:t_stat} is fundamental to proving \cref{th:loco_t}.

\begin{proof}
Proving \cref{th:t_stat} begins with the definition for the $t$-statistic:
\begin{equation}
\tag{\ref{eq:t_stat_og}}
    t_i = \frac{\hat{\beta}_i}{SE(\hat{\beta}_i)} = \frac{\hat{\beta}_i}{\sqrt{\operatorname{Var}(\epsilon)(\mathbf{X^{T}X})^{-1}_{ii}}}.
\end{equation}

The main bottleneck here is $\mathbf{X^{T}X}^{-1}_{ii}$. To address this, consider the definition of a covariance matrix is $cov(\mathbf{X}) = \frac{\mathbf{X^{T}X}}{n-1}$. Additionally, recall from \cref{eq:cov_x} (found in \cref{ss2:latent}) that $cov(\mathbf{X}) = \mathbf{AA^T} = (2\Delta + (p-2)\Delta^2)\mathbf{J}+(1-\Delta)^2\mathbf{I}$. If we let $\alpha = \frac{2\Delta + (p-2)\Delta^2}{(1-\Delta)^2}$, then we can combine definitions with a simplification of \cref{eq:cov_x} to show that
\begin{equation}
\label{eq:xtx}
    \mathbf{X^T X} = (n-1)\mathbf{AA^T} = (n-1)(1-\Delta)^2 [\mathbf{I} + \alpha\mathbf{J}].
\end{equation}

To solve for $\mathbf{X^{T}X}^{-1}_{ii}$, we use the Sherman–Morrison formula \citep{sherman_morrison} which states:
\begin{equation}
\label{eq:sherman_morrison}
        (\mathbf{M} + uv^T)^{-1}= \mathbf{M}^{-1} -  \frac{\mathbf{M}^{-1}uv^T\mathbf{M}^{-1}}{1 + v^T \mathbf{M}^{-1} u}.
\end{equation}

Since $\alpha \mathbf{J}$ in \cref{eq:xtx} is rank-one, we can express it as a multiplication of vectors filled with 1s: $\alpha \mathbf{1}\mathbf{1}^T$. Now, let $\mathbf{I} = \mathbf{M} = \mathbf{M}^{-1}$ and let $u = v =\sqrt{\alpha \mathbf{1}}$ be vectors of length $p$ with all entries $\sqrt{\alpha}$. Then we can find $(\mathbf{I} + \alpha\mathbf{J})^{-1}$ via the Sherman Morrison formula with some variable substitution and simplification which yields
\begin{equation}
\label{eq:inverse_mat}
    (\mathbf{I} + \alpha\mathbf{J})^{-1}= (\mathbf{I} + uv^T)^{-1}= \mathbf{I} -  \frac{uv^T}{1 + v^T u} = \mathbf{I} -  \frac{\alpha \mathbf{J}}{1 + \alpha p}.
\end{equation}

Therefore,
\begin{equation}
\label{eq:xtxi}
    \begin{split}
        \mathbf{X^{T}X}^{-1} &= ((n-1)(1-\Delta)^2 [\mathbf{I} + \alpha\mathbf{J}])^{-1} \\
        &= \frac{1}{(n-1)(1-\Delta)^2}\left[\mathbf{I} -  \frac{\alpha \mathbf{J}}{1 + \alpha p}\right].
    \end{split}
\end{equation}

Substituting $\mathbf{I}_{ii} = \mathbf{J}_{ii} = 1$ into \cref{eq:xtxi} and simplifying yields:
\begin{equation}
\label{eq:xtxi_simp1}
\begin{split}
     (\mathbf{X^{T}X})^{-1}_{ii} &= \frac{1}{(n-1)(1-\Delta)^2}\left[1 -  \frac{\alpha}{1 + \alpha p}\right]\\
     &= \frac{1}{(n-1)(1-\Delta)^2}\left[\frac{1 + \alpha(p-1)}{1 + \alpha p}\right].
\end{split}
\end{equation}

Substituting $\alpha$ into the latter expression of \cref{eq:xtxi_simp1} and simplifying yields:
\begin{equation}
\label{eq:xtxi_simp2}
\begin{split}
     \left[\frac{1 + (p-1)\alpha}{1 + p\alpha }\right] &= \left[\frac{1 + (p-1)\frac{2\Delta + (p-2)\Delta^2}{(1-\Delta)^2}}{1 + p\frac{2\Delta + (p-2)\Delta^2}{(1-\Delta)^2}}\right]\\
     &= \left[\frac{\frac{(1-\Delta)^2 + 2(p-1)\Delta + (p-1)(p-2)\Delta^2}{(1-\Delta)^2}}{\frac{(1-\Delta)^2 + 2p\Delta + p(p-2)\Delta^2}{(1-\Delta)^2}}\right]\\
     &= \left[\frac{(1-\Delta)^2 + 2(p-1)\Delta + (p-1)(p-2)\Delta^2}{(1-\Delta)^2 + 2p\Delta + p(p-2)\Delta^2}\right]\\
     &= \left[\frac{1-2\Delta +\Delta^2 + 2p\Delta - 2\Delta + p^2\Delta^2 -3p\Delta^2 + 2\Delta^2}{1-2\Delta +\Delta^2 + 2p\Delta + p^2\Delta^2 -2p\Delta^2}\right]\\
     &= \left[\frac{1 + 2p\Delta - 4\Delta + p^2\Delta^2 -3p\Delta^2 + 3\Delta^2}{1 + 2p\Delta - 2\Delta + p^2\Delta^2 -2p\Delta^2 +\Delta^2}\right]\\
     &= \left[\frac{1 + 2p\Delta - 4\Delta + p^2\Delta^2 -4p\Delta^2  + 4\Delta^2 +(p-1)\Delta^2}{1 + 2(p-1)\Delta + (p-1)^2\Delta^2}\right]\\
     &= \left[\frac{1 + 2(p-2)\Delta + (p-2)\Delta^2 +(p-1)\Delta^2}{(1 + (p-1)\Delta)^2}\right]\\
     &= \left[\frac{(1 + (p-2)\Delta)^2 + (p-1)\Delta^2}{(1 + (p-1)\Delta)^2}\right].
\end{split}
\end{equation}

Therefore,
\begin{equation}
\label{eq:xtxi_simp3}
\begin{split}
    (\mathbf{X^{T}X})^{-1}_{ii} &= \frac{1}{(n-1)(1-\Delta)^2}\left[\frac{1 + \alpha(p-1)}{1 + \alpha p}\right]\\
    &= \frac{(1 + (p-2)\Delta)^2 + (p-1)\Delta^2}{(n-1)(1-\Delta)^2(1 +(p-1)\Delta)^2}.
\end{split}
\end{equation}

Substituting \cref{eq:xtxi_simp3} into \cref{eq:t_stat_og} and assuming $\beta=\hat{\beta}$ yields
\begin{equation}
\tag{\ref{eq:t_stat}}
    t_i = \beta_i(1-\Delta) \sqrt{\frac{(n-1)(1+(p-1)\Delta)^2}{\operatorname{Var}(\epsilon)[(1 + (p-2)\Delta)^2 + (p-1)\Delta^2]}}. \qed
\end{equation} 
\end{proof}

\section{Empirical Relationship between LOCO and \texorpdfstring{$t$}{t}-statistic}
\label{app:t_results}

\begin{figure}[htbp]
\begin{center}
\centerline{\includegraphics[width=0.75\columnwidth]{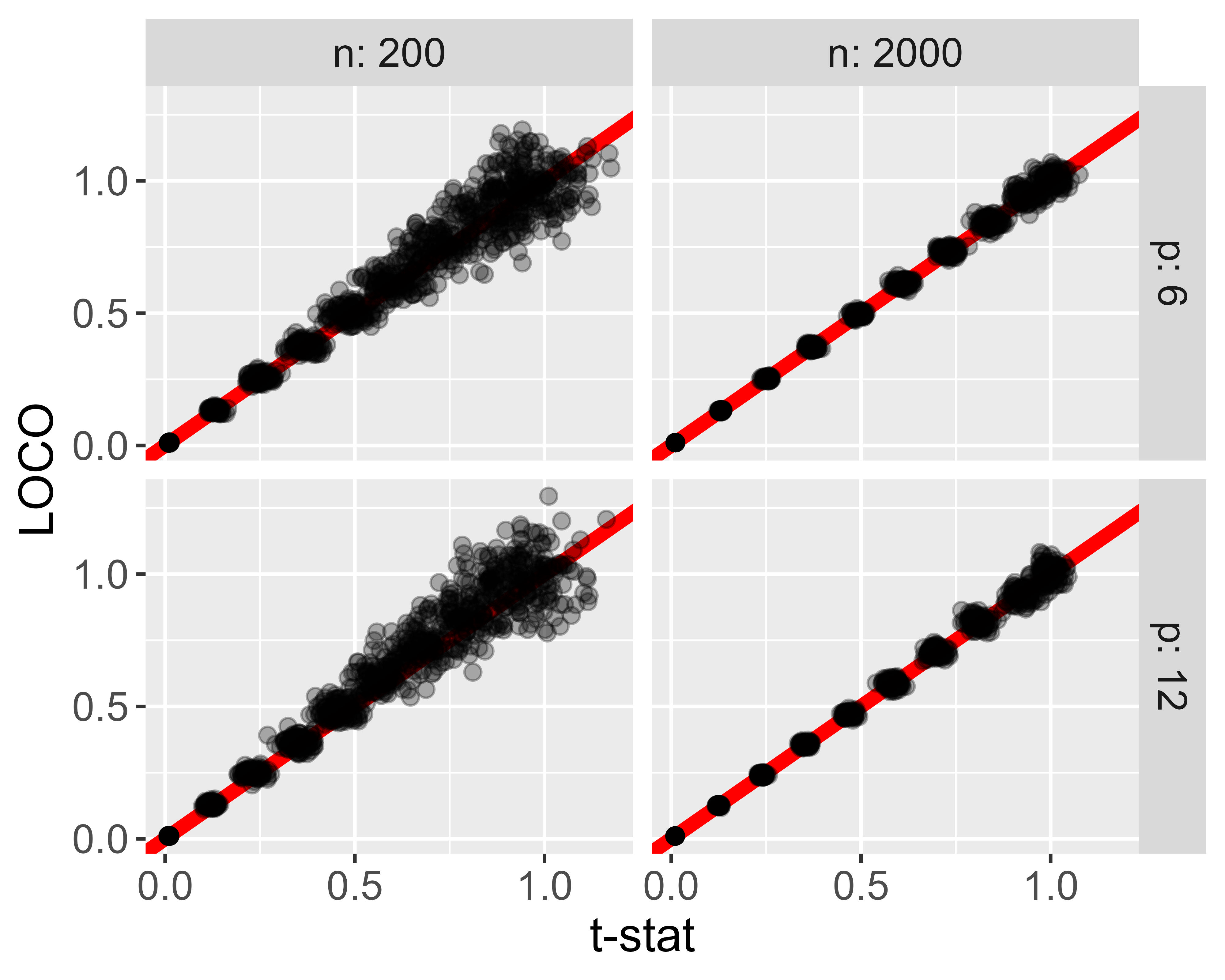}}
\caption{Parity scatter plots for assessing the relationship between empirical $t$-statistics and empirical LOCO values for data defined in \cref{ss2:data}. Red lines denote a slope of 1. Plots are faceted by $p$ and $n$ dimensions. We used \cref{eq:loco_t} to ensure a slope of 1 by transforming the $t$-statistics. We also applied a small sample bias correction from \cref{s3:ss_bias} to LOCO. No small sample bias correction was applied to the $t$-statistics since this is already accounted for in their calculation.}
\label{fig:t_vs_loco}
\end{center}
\end{figure}

In similar fashion to \cref{ss3:lin_imp}, we offer \cref{fig:t_vs_loco} for comparing empirical $t$-statistics to empirical LOCO metrics and validating our derivation of \cref{eq:loco_t}. Since both metrics are empirical, we employ a parity scatter plots to display this comparison. We apply a small sample bias correction from \cref{s3:ss_bias} to LOCO. Conversely, no small sample bias correction is applied to the $t$-statistics since this is already accounted for in their calculation. Noticeable clusters are monotonically driven by $\Delta$, where empirical values near 1 correspond to low $\Delta$ and values near 0 correspond to high $\Delta$.

Since LOCO values naturally ranged from 0 to 1, we use and reverse \cref{eq:loco_t} to multiply the $t$-statistics by $\sqrt{\frac{\operatorname{Var}(\epsilon)}{n-1}}$, with the understanding that the scatter plots should have a slope of 1 if \cref{eq:loco_t} is correct. 

From the comparisons displayed in \cref{fig:t_vs_loco}, \cref{eq:loco_t} appears to be correct, as the resulting points visually fit a line with a slope of 1 extremely well. This offers empirical evidence supporting the work in \cref{ss2:loco} and \cref{app:t_stat}.

\section*{Acknowledgments}
The authors thank Peter Crooks for reviewing the validity and relevance of our derivation relating LOCO to the $t$-statistic and for his suggestions of how to highlight the legitimacy of our $c$ definition through empirical demonstration. The authors also thank Kevin Moon for recommending promising journals and venues for this work and Adam Robertson for providing helpful feedback during the development of the LOCO derivation.

The authors also appreciate OpenAI’s ChatGPT (GPT-5) for assistance with literature review, formula cross-checking, and language editing for improving the readability of the manuscript. While its help was greatly appreciated, all research contributions, analyses, and conclusions are solely the responsibility of the authors.


The code used to generate the data and figures in this article can be found at:
\url{https://github.com/KelvynBladen/theoryVarImportance}

\bibliographystyle{references/siamplain}
\bibliography{main}

\end{document}